\documentclass[11pt]{amsart}
\usepackage{amssymb}
\usepackage{amsaddr}
\usepackage[margin=1in]{geometry}
\usepackage[colorlinks]{hyperref}
\usepackage{graphicx}

\newtheorem{theorem}{Theorem}[section]

\newtheorem{prop}[theorem]{Proposition}
\newtheorem{cor}[theorem]{Corollary}

\theoremstyle{definition}

\newtheorem{example}[theorem]{Example}

\theoremstyle{remark}
\newtheorem{remark}[theorem]{Remark}

\numberwithin{equation}{section}

\begin{document}

\newcommand{\spacing}[1]{\renewcommand{\baselinestretch}{#1}\large\normalsize}
\spacing{1.14}
\title[Ricci Solitons of Left-Invariant Metrics]{Ricci Solitons of Left-Invariant Metrics on Lie Groups with Derived Algebra of Dimension at Most Two}

\author {Hamid Reza Salimi Moghaddam}

\address{Department of Pure Mathematics, Faculty of  Mathematics and Statistics, University of Isfahan, Isfahan, 81746-73441-Iran.
Scopus Author ID: 26534920800 \\ ORCID Id:0000-0001-6112-4259}
\email{salimi.moghaddam@gmail.com and hr.salimi@sci.ui.ac.ir}
\keywords{left-invariant Riemannian metric, Ricci curvatures, derived algebra, Ricci soliton.\\
\hspace*{0.4cm} AMS 2020 Mathematics Subject Classification: 53C30, 53C21, 22E60.}

\date{\today}

\begin{abstract}
In this article, we investigate when a left-invariant Riemannian metric on a Lie group is a Ricci soliton,
under the assumption that the derived algebra has dimension at most two. We establish computable necessary and sufficient
conditions for a given left-invariant Riemannian metric to be a Ricci soliton. As applications, we obtain several examples
of Ricci nilsolitons and apply our results to indecomposable Lie groups of dimension at least five with two-dimensional derived algebra.\\
\end{abstract}

\maketitle

\section{\textbf{Introduction}}\label{Introduction}
In 1982, Hamilton initiated the study of the Ricci flow in the context of three-dimensional Riemannian manifolds. In \cite{Hamilton}, he demonstrated that when the Ricci curvature of such manifolds is strictly positive, they also admit a Riemannian metric of constant positive curvature. Hamilton's seminal work focused on the unnormalized Ricci flow equation
\begin{equation}\label{Ricci flow equation}
    \frac{\partial}{\partial t}g_t=-2\textsf{ric}_{g(t)},
\end{equation}
where $g_t$ represents a one-parameter family of Riemannian metrics. A Riemannian metric $g$ is called a Ricci soliton if there exists a real number $c$ and a complete vector field $X$ such that:
\begin{equation}\label{Main Ricci soliton equation}
    \textsf{ric}=cg+\textsf{L}_Xg.
\end{equation}
On the other hand, if for a one-parameter family of diffeomorphisms $\phi_t$, and a family of scalars $\tau(t)$, on an interval $(t_1,t_2)$ containing zero, $(M,g_t=\tau(t)\phi^{\ast}_tg_0)$ be a solution of the evolution equation \eqref{Ricci flow equation}, it is termed a self-similar solution. Using the correspondence between Ricci solitons and self-similar solutions (Lemma 2.4 of \cite{Chow-Knopf}), they are regarded as equivalent. This family contains the class of Einstein manifolds. In equation \eqref{Main Ricci soliton equation}, if $c>0$, $c=0$, and $c<0$, the Ricci soliton is referred to as shrinking, steady, and expanding, respectively.

As in other geometric settings, an interesting case arises when investigating Lie groups with left-invariant metrics. Assuming $G$ is a Lie group with Lie algebra $\frak{g}$, a left-invariant Riemannian metric $g$ on $G$ is called an algebraic Ricci soliton if there exists a real number $c$ and a derivation $D\in\textsf{Der}(\frak{g})$ such that:
\begin{equation}\label{Algebraic Ricci soliton equation}
    \textsf{Ric}=c.\operatorname{Id}+D,
\end{equation}
where $\textsf{Ric}$ represents the Ricci transformation of $g$.

In 2001, Lauret showed that Ricci solitons on nilpotent Lie groups are semi-algebraic (see \cite{Lauret1}). In 2014, Jablonski extended Lauret's results to homogeneous spaces, demonstrating that all homogeneous Ricci solitons are algebraic (see \cite{Jablonski2}). For a complete survey on this field, see \cite{Jablonski3}.

We recall that a Ricci nilsoliton $(N,\langle \cdot , \cdot \rangle)$ is a nilpotent Lie group $N$ equipped with a left-invariant Riemannian metric satisfying equation \eqref{Algebraic Ricci soliton equation}. Let $\langle \cdot , \cdot \rangle$ be a left-invariant Riemannian metric on a nilpotent Lie group $N$ with Lie algebra $\frak{n}$. In this case, a metric solvable Lie algebra $(\frak{s}=\frak{a}\oplus\frak{n},\langle \cdot , \cdot \rangle')$ is called a metric solvable extension of $(N,\langle \cdot , \cdot \rangle)$ if the following relations hold:
\begin{equation*}
    \frak{s}^1=[\frak{s},\frak{s}]=\frak{n}=\frak{a}^\perp, \ \ \ [\cdot,\cdot]_\frak{s}|_{\frak{n}\times\frak{n}}=[\cdot,\cdot]_\frak{n}, \ \ \ \langle \cdot , \cdot \rangle'|_{\frak{n}\times\frak{n}}=\langle \cdot , \cdot \rangle,
\end{equation*}
where we have denoted the Lie brackets of $\frak{s}$ and $\frak{n}$ by $[\cdot,\cdot]_\frak{s}$ and $[\cdot,\cdot]_\frak{n}$, respectively. If $\frak{a}=\frak{n}^\perp$ is abelian, then such an extension is called standard. In \cite{Lauret1}, Lauret showed that a nilmanifold $(N,\langle \cdot , \cdot \rangle)$ is a Ricci nilsoliton if and only if it admits an Einstein standard solvable extension. On the other hand, in \cite{Lauret2}, he proved that any Einstein solvmanifold is given by such a standard extension. Thus, the classification of Einstein solvmanifolds is equivalent to the classification of Ricci nilsolitons.
In \cite{Nikolayevsky2}, Nikolayevsky gave a rather complicated criterion for a two-step nilpotent Lie group with two-dimensional derived algebra to be a Ricci nilsoliton. Recently, Gordon and Jablonski proved that left-invariant Ricci solitons on solvable Lie groups, up to scaling and isometry, are unique (see \cite{Gordon-Jablonski2024}).

In the present paper, rather than addressing the existence of a left-invariant Riemannian metric that is a Ricci soliton, we consider the following question: given a left-invariant Riemannian metric on a Lie group, when is the metric a Ricci soliton? We study the general case of Lie groups with two-dimensional derived algebra, while in Corollary \ref{corollary two-step nilpotent two-dim commutator} we give a necessary and sufficient condition for a two-step nilpotent Lie group with two-dimensional derived algebra to be a Ricci nilsoliton. This criterion can be easily checked using the structure constants of the Lie algebra with respect to an orthonormal basis (see Examples \ref{example 1 two-step nilpotent} and \ref{example 2 two-step nilpotent}).

We first study left-invariant Riemannian metrics on Lie groups with one- and two-dimensional derived algebras and establish necessary and sufficient conditions for such a metric to be a Ricci soliton. Then, we give several examples of Ricci nilsolitons. Finally, we apply these criteria to indecomposable Lie groups of dimension at least five with two-dimensional derived algebra and show that, except for one six-dimensional case, the corresponding left-invariant Riemannian metrics are not Ricci solitons.

\textbf{Notations:} In this work, we use $\textsf{Lie(n,m)}$, ($m=1, 2$), as the class of all $n$-dimensional real Lie groups with $m$-dimensional derived algebra. For a real Lie group $G$ with Lie algebra $\frak{g}$, we denote the derived algebra of $\frak{g}$ by $\frak{g}^1=[\frak{g},\frak{g}]$, and by $\textsf{Der}(\frak{g})$ the Lie algebra of derivations of $\frak{g}$. Also, we use the same notation $\langle \cdot , \cdot \rangle$ for a left-invariant Riemannian metric on a Lie group $G$ and the induced inner product on the Lie algebra $\frak{g}$. For a linear map $f$ we use the notation $f^m$ for $\overbrace{f\circ\dots\circ f}^{m\rm\ times}$.


\section{\textbf{Left-invariant Riemannian metrics on Lie groups with one-dimensional derived algebras}}\label{1-dimensional}
In this section, we study the geometry of left-invariant Riemannian metrics on Lie groups $G$, where $G\in\textsf{Lie(n,1)}$.
Up to local isomorphism, the set $\textsf{Lie(n,1)}$ has exactly $[\frac{n+1}{2}]$ different elements. In fact, a Lie group $G\in\textsf{Lie(n,1)}$ is of the form $G=H\times K$, where $H$ is a $(2m+1)$-dimensional Heisenberg group, $m=1,2,\dots, [\frac{n-1}{2}]$, or the affine group $A(1)$, and $K$ is an abelian group (see \cite{Kaiser}). Kaiser studied the Ricci curvatures of left-invariant Riemannian metrics on Lie groups belonging to $\textsf{Lie(n,1)}$ in \cite{Kaiser}. For an arbitrary Lie group $G\in\textsf{Lie(n,1)}$ equipped with an arbitrary left-invariant Riemannian metric $\langle \cdot , \cdot \rangle$, we review the Levi-Civita connection, sectional curvature, and Ricci curvature, and then give a necessary and sufficient condition for the metric to be a Ricci soliton.

Let $G\in\textsf{Lie(n,1)}$ and $\langle \cdot , \cdot \rangle$ be a left-invariant Riemannian metric on $G$. Suppose that $\textsf{e}\in\frak{g}^1$ is a unit vector and $\Gamma\subset\frak{g}$ is the hyperplane orthogonal to the vector $\textsf{e}$, with respect to $\langle \cdot , \cdot \rangle$. So there exists a linear map $\phi:\Gamma\to\Bbb{R}$ such that for every $u\in\Gamma$ we have $[u,\textsf{e}]=\phi(u)\textsf{e}$.
Hence there is a unique element $\textsf{a}\in\Gamma$ such that $\phi(u)=\langle\textsf{a},u\rangle$, for all $u\in\Gamma$. In fact for any $u\in\Gamma$ we have $[u,\textsf{e}]=\langle\textsf{a},u\rangle\textsf{e}$. On the other hand, there is a skew-symmetric bilinear form $B:\Gamma\times\Gamma\to\Bbb{R}$ such that $[u,v]=B(u,v)\textsf{e}$, for all $u, v\in\Gamma$. So there exists a unique skew-adjoint (with respect to $\langle \cdot , \cdot \rangle$) linear map $f:\Gamma\to\Gamma$ such that for all $u,v\in\Gamma$ we have $B(u,v)=\langle f(u) , v \rangle$ or equivalently $[u,v]=\langle f(u) , v \rangle\textsf{e}$.

It is easy to see that for the Levi-Civita connection of the Riemannian manifold $(G,\langle \cdot , \cdot \rangle)$ we have (see \cite{Asgari-Salimi} and \cite{Kaiser}):
\begin{equation*}
    \begin{array}{cc}
  \nabla_\textsf{e}\textsf{e}=\textsf{a}, & \nabla_\textsf{e}u=-\frac{1}{2}f(u)-\langle u , \textsf{a} \rangle\textsf{e}, \\
  \nabla_u\textsf{e}=-\frac{1}{2}f(u), & \nabla_uv=\frac{1}{2}\langle f(u) , v \rangle\textsf{e},
    \end{array}
\end{equation*}
for all $u, v\in\Gamma$. Also, for the sectional curvature, we have (see \cite{Asgari-Salimi}):
\begin{equation*}
    K(u,\textsf{e})=\frac{1}{4}\|f(u)\|^2-\langle \textsf{a} , u \rangle^2, \ \ \ \ \ K(u,v)=-\frac{3}{4}\langle f(u) , v \rangle^2,
\end{equation*}
where $\{u,v\}$  is an orthonormal set of elements of $\Gamma$.

The Ricci transformation of $(G,\langle \cdot , \cdot \rangle)$ is of the following form (see \cite{Kaiser}):
\begin{eqnarray}
  \textsf{Ric}(\textsf{e}) &=& -f(\textsf{a})-(\frac{1}{4}\operatorname{tr}(f^2)+\|\textsf{a}\|^2)\textsf{e}, \label{Ricci One-dim-1} \\
  \textsf{Ric}(u) &=& \frac{1}{2}f^2(u)-\langle u , \textsf{a} \rangle\textsf{a}+\langle f(u) , \textsf{a} \rangle\textsf{e}.\label{Ricci One-dim-2}
\end{eqnarray}
After the above overview of the geometry of $(G,\langle . , .\rangle)$ for $G\in\textsf{Lie(n,1)}$, we present the main theorem of this section.

\begin{theorem}\label{Classi One-dim commu Ricci solitons}
A Lie group $G\in\textsf{Lie(n,1)}$ equipped with an arbitrary left-invariant Riemannian metric $\langle\cdot,\cdot\rangle$ is a Ricci soliton if and only if one of the following two cases holds:
\begin{description}
\item[i] $G$ is unimodular, i.e. $\textsf{a}=0$, and there exists $c\in\mathbb{R}$ such that

$$
f^3=\left(c-\frac14\operatorname{tr}(f^2)\right)f.
$$

In this case, $G$ is a two-step nilpotent Lie group and $(G,\langle\cdot,\cdot\rangle)$ is a Ricci nilsoliton. Since $\mathfrak{g}^1\neq0$, one has $f\neq0$, and hence $c<0$.

\item[ii] $G$ is non-unimodular, $0\neq\textsf{a}\in\ker(f)$,

$$
c=-\|\textsf{a}\|^2,
$$

and

$$
f^3=-\left(2\|\textsf{a}\|^2+\frac14\operatorname{tr}(f^2)\right)f.
$$

In this case, the Ricci soliton is expanding.
\end{description}
\end{theorem}
\begin{proof}
We know that $(G,\langle\cdot,\cdot\rangle)$ is a Ricci soliton if and only if there exist $c\in\mathbb{R}$ and $D\in\operatorname{Der}(\mathfrak g)$ such that:

$$
\textsf{Ric}(x)=cx+D(x),\qquad x\in\mathfrak g.
$$

Hence, by \eqref{Ricci One-dim-1} and \eqref{Ricci One-dim-2},

$$
D(\textsf e)
=
-f(\textsf a)
-\left(\frac14\operatorname{tr}(f^2)+\|\textsf a\|^2+c\right)\textsf e,
$$

and

$$
D(u)
=
\frac12f^2(u)-\langle u,\textsf a\rangle\textsf a
+\langle f(u),\textsf a\rangle\textsf e-cu.
$$

Since $D$ is a derivation, applying it to $[u,\textsf e]=\langle\textsf a,u\rangle\textsf e$ gives
\begin{equation}\label{[D(u),e] One-dim-4}
   \langle \textsf{a} , \frac{3}{2}f^2(u)-\|\textsf{a}\|^2u-cu \rangle = 0,
\end{equation}
and
\begin{equation}\label{[D(u),e] One-dim-5}
    \langle \textsf{a} , u \rangle f(\textsf{a})=0.
\end{equation}
for every $u\in\Gamma$, and hence, if $\textsf a\neq0$,

$$
f(\textsf a)=0.
$$

In the same case, equation \eqref{[D(u),e] One-dim-4}, with $u=\textsf a$, gives

$$
c=-\|\textsf a\|^2.
$$

Substituting this value of $c$ into \eqref{[D(u),e] One-dim-4}, we see that the latter equation is then automatically satisfied for every $u\in\Gamma$.

Now suppose first that $\textsf a=0$. Then $G$ is unimodular and

$$
D(\textsf e)
=
-\left(\frac14\operatorname{tr}(f^2)+c\right)\textsf e,
\qquad
D(u)=\frac12f^2(u)-cu.
$$

Applying the derivation identity to $[u,v]=\langle f(u),v\rangle\textsf e$, we obtain

$$
-\left(\frac14\operatorname{tr}(f^2)+c\right)
\langle f(u),v\rangle
=
\langle f^3(u)-2cf(u),v\rangle.
$$

Therefore,

\begin{equation}\label{a=0}
    f^3=
\left(c-\frac14\operatorname{tr}(f^2)\right)f.
\end{equation}

Conversely, if $\textsf a=0$ and this identity holds for some $c$, the above formulas define a derivation $D$, and hence the metric is a Ricci soliton. Since $f\neq0$ and $f$ is skew-adjoint, we have $\operatorname{tr}(f^2)<0$ and $\operatorname{tr}(f^4)>0$. So the equation \eqref{a=0} implies that $\operatorname{tr}(f^4)=(c-\frac{1}{4}\operatorname{tr}(f^2))\operatorname{tr}(f^2)$. Hence the corresponding constant $c$ is negative which shows that the soliton is expanding and is a Ricci nilsoliton.

Finally, suppose that $\textsf a\neq0$. Applying the derivation identity to $[u,v]$ and using $c=-|\textsf a|^2$, we obtain

$$
f^3
=
-\left(
2\|\textsf a\|^2+\frac14\operatorname{tr}(f^2)
\right)f.
$$

Conversely, this identity together with $f(\textsf a)=0$ and $c=-|\textsf a|^2$ ensures that the above map $D$ is a derivation. Hence the metric is an expanding Ricci soliton.
\end{proof}

\begin{remark}
Jablonski and Nikolayevsky, using different techniques, proved that a nilpotent Lie group $G=H\times K$ is a Ricci nilsoliton if and only if both $H$ and $K$ are Ricci nilsolitons (see \cite{Jablonski1} and \cite{Nikolayevsky1}). Therefore, to characterize the Ricci nilsoliton structures on a two-step nilpotent Lie group $G\in\textsf{Lie(n,1)}$, it is sufficient to classify Ricci nilsoliton structures on the Lie groups $\operatorname{Heis}_{2m+1}$ and $\Bbb{R}^{n-(2m+1)}$.
\end{remark}

\section{\textbf{Left-invariant Riemannian metrics on Lie groups with two-dimensional derived algebras}}\label{2-dimensional}
This section considers the Lie groups $G\in\textsf{Lie(n,2)}$, ($n\geq 3$). Since $G$ is solvable, $\frak{g}^1$ is a two-dimensional nilpotent Lie algebra; hence the derived algebra is abelian.
Therefore, $\frak{g}$ is a 2-step solvable Lie algebra (for more details see \cite{Le-Nguyen-Vo}). As in the previous section, let $\langle \cdot , \cdot \rangle$ be an arbitrary left-invariant Riemannian metric on $G$.
Assume that $\{\textsf{e}_1,\textsf{e}_2\}$ is an orthonormal basis for $\frak{g}^1$ and $\Gamma$ is the $(n-2)$-dimensional subspace of $\frak{g}$ which is orthogonal to $\frak{g}^1$, with respect to the inner product $\langle \cdot , \cdot \rangle$. Then there exist four real linear maps $\phi_1, \phi_2, \psi_1, \psi_2:\Gamma\to\Bbb{R}$ such that for any $u\in\Gamma$ we have:
\begin{eqnarray*}
  {[}u,\textsf{e}_1{]} &=& \phi_1(u)\textsf{e}_1+\phi_2(u)\textsf{e}_2, \\
  {[}u,\textsf{e}_2{]} &=& \psi_1(u)\textsf{e}_1+\psi_2(u)\textsf{e}_2.
\end{eqnarray*}
So there exist unique vectors $\textsf{a}_1, \textsf{a}_2, \textsf{b}_1, \textsf{b}_2 \in\Gamma$ such that for any $u\in\Gamma$ we have:
\begin{equation*}
    \phi_1(u)=\langle\textsf{a}_1,u\rangle, \phi_2(u)=\langle\textsf{a}_2,u\rangle, \psi_1(u)=\langle\textsf{b}_1,u\rangle, \psi_2(u)=\langle\textsf{b}_2,u\rangle.
\end{equation*}
Applying the Jacobi identity to the triple $(\textsf{a}_2, \textsf{b}_1, \textsf{e}_1)$ (or equivalently to $(\textsf{a}_2, \textsf{b}_1, \textsf{e}_2)$), along with the Cauchy-Schwarz inequality, leads to the conclusion that $\textsf{b}_1$ is a scalar multiple of $\textsf{a}_2$. However, we will not use this fact to establish the symmetry of the equations.

As in the previous section, there are two skew-symmetric bilinear forms $B_1, B_2:\Gamma\times\Gamma\to\Bbb{R}$ such that for all $u, v\in\Gamma$ we have $[u,v]=B_1(u,v)\textsf{e}_1+B_2(u,v)\textsf{e}_2$.
So there are two unique skew-adjoint (with respect to $\langle \cdot , \cdot \rangle$) linear maps $f_1, f_2:\Gamma\to\Gamma$ such that $B_1(u,v)=\langle f_1(u) , v \rangle$ and $B_2(u,v)=\langle f_2(u) , v \rangle$,
for all $u,v\in\Gamma$. Note that for the two endomorphisms $f_1$ and $f_2$ we have, $\operatorname{tr}(f_1\circ f_2)=\operatorname{tr}(f_2\circ f_1)$.

Using the formula $\nabla_xy=\frac{1}{2}(\textsf{ad}_xy-\textsf{ad}_x^\ast y-\textsf{ad}_y^\ast x)$ for left-invariant Riemannian metrics, together with the sectional curvature formula, we obtain the following proposition.

\begin{prop}\label{Levi-Civita and Sectional}
The Levi-Civita connection of any Lie group $G\in\textsf{Lie(n,2)}$ equipped with a left-invariant Riemannian metric $\langle \cdot , \cdot \rangle$,  is given in Table \ref{table 1}, where $u$ and $v$ are arbitrary elements of $\Gamma$.
\begin{table}[h]
\caption{Levi-Civita connection of a Lie group $G\in\textsf{Lie(n,2)}$ equipped with a left-invariant Riemannian metric}
\label{table 1}
\centering
  \begin{tabular}{|c|c|c|c|}
    \hline
    $\nabla$ & $\textsf{e}_1$ & $\textsf{e}_2$ & $u$ \\
    \hline
    $\textsf{e}_1$ & $\textsf{a}_1$ & $\frac{1}{2}(\textsf{a}_2+\textsf{b}_1)$ & $-\langle\textsf{a}_1,u\rangle\textsf{e}_1-\frac{1}{2}(\langle\textsf{b}_1+\textsf{a}_2,u\rangle\textsf{e}_2+f_1(u))$ \\
    \hline
    $\textsf{e}_2$ & $\frac{1}{2}(\textsf{a}_2+\textsf{b}_1)$ & $\textsf{b}_2$ & $-\langle\textsf{b}_2,u\rangle\textsf{e}_2-\frac{1}{2}(\langle\textsf{b}_1+\textsf{a}_2,u\rangle\textsf{e}_1+f_2(u))$ \\
    \hline
    $v$ & $\frac{1}{2}(\langle\textsf{a}_2-\textsf{b}_1,v\rangle\textsf{e}_2-f_1(v))$ & $\frac{1}{2}(\langle\textsf{b}_1-\textsf{a}_2,v\rangle\textsf{e}_1-f_2(v))$ & $\frac{1}{2}(\langle f_1(v),u\rangle\textsf{e}_1+\langle f_2(v),u\rangle\textsf{e}_2)$ \\
    \hline
\end{tabular}
\end{table}

For the sectional curvature we have:
\begin{eqnarray*}
    K(\textsf{e}_1,\textsf{e}_2) &=& \frac{1}{4}\|\textsf{a}_2+\textsf{b}_1\|^2-\langle\textsf{a}_1,\textsf{b}_2\rangle,\\
    K(\textsf{e}_1,u) &=& \frac{1}{4}\Big{(}\langle\textsf{b}_1,u\rangle^2-3\langle\textsf{a}_2,u\rangle^2+\|f_1(u)\|^2\Big{)}-\langle\textsf{a}_1,u\rangle^2-\frac{1}{2}\langle\textsf{a}_2,u\rangle\langle\textsf{b}_1,u\rangle, \\
    K(\textsf{e}_2,u) &=& \frac{1}{4}\Big{(}\langle\textsf{a}_2,u\rangle^2-3\langle\textsf{b}_1,u\rangle^2+\|f_2(u)\|^2\Big{)}-\langle\textsf{b}_2,u\rangle^2-\frac{1}{2}\langle\textsf{a}_2,u\rangle\langle\textsf{b}_1,u\rangle, \\
    K(u,v) &=& -\frac{3}{4}\Big{(}\langle u,f_1(v)\rangle^2+\langle u,f_2(v)\rangle^2\Big{)},
\end{eqnarray*}
where $\{u, v\}$ is an orthonormal set in $\Gamma$.
\end{prop}

The main goal of this section is to give a computable necessary and sufficient condition for $(G,\langle \cdot , \cdot \rangle)$ to be a Ricci soliton, where $G\in\textsf{Lie(n,2)}$ and $n\geq 3$.
To this end, we define the following quantities. These quantities will appear in the next theorem.
\begin{eqnarray}
  A_1 &:=& \frac{1}{2}(\|\textsf{b}_1\|^2-\|\textsf{a}_2\|^2)-\frac{1}{4}\operatorname{tr}(f^2_1)-\|\textsf{a}_1\|^2-\langle\textsf{a}_1,\textsf{b}_2\rangle, \\
  B_1 &:=& -(\langle\textsf{a}_1,\textsf{b}_1\rangle+\langle\textsf{a}_2,\textsf{b}_2\rangle+\frac{1}{4}\operatorname{tr}(f_1\circ f_2)), \\
  E_1 &:=& -f_1(\textsf{a}_1+\frac{1}{2}\textsf{b}_2)-\frac{1}{2}f_2(\textsf{a}_2), \\
  A_2 &:=& B_1 \\
  B_2 &:=& \frac{1}{2}(\|\textsf{a}_2\|^2-\|\textsf{b}_1\|^2)-\frac{1}{4}\operatorname{tr}(f^2_2)-\|\textsf{b}_2\|^2-\langle\textsf{a}_1,\textsf{b}_2\rangle, \\
  E_2 &:=& -f_2(\textsf{b}_2+\frac{1}{2}\textsf{a}_1)-\frac{1}{2}f_1(\textsf{b}_1).
\end{eqnarray}
Also, for any $u\in\Gamma$ we define:
\begin{eqnarray}
  A(u) &:=& -\langle u,f_1(\textsf{a}_1+\frac{1}{2}\textsf{b}_2)+\frac{1}{2}f_2(\textsf{a}_2)\rangle, \\
  B(u) &:=& -\langle u,f_2(\textsf{b}_2+\frac{1}{2}\textsf{a}_1)+\frac{1}{2}f_1(\textsf{b}_1)\rangle, \\
  E(u) &:=& -\langle u,\textsf{a}_1\rangle\textsf{a}_1-\langle u,\textsf{b}_2\rangle\textsf{b}_2-\frac{1}{2}\langle u,\textsf{b}_1+\textsf{a}_2\rangle(\textsf{b}_1+\textsf{a}_2)+\frac{1}{2}(f_1^2+f_2^2)(u).
\end{eqnarray}

We now have all the tools needed to state the main theorem of this section.
\begin{theorem}\label{Classi Two-dim commu Ricci solitons}
A Lie group $G\in\textsf{Lie(n,2)}$ equipped with an arbitrary left-invariant Riemannian metric $\langle . , .\rangle$ is a Ricci soliton if and only if there exists $c\in\Bbb{R}$ such that the following equations hold for any $u, v\in\Gamma$.
\begin{eqnarray}
 && \langle \textsf{b}_1, E_1\rangle-\langle \textsf{a}_1, E_2\rangle = 0, \label{eq1}\\
 && \langle \textsf{b}_2, E_1\rangle-\langle \textsf{a}_2, E_2\rangle = 0, \label{eq2}\\
 && \langle \textsf{a}_1, u\rangle E_1+\langle \textsf{a}_2, u\rangle E_2 = 0, \label{eq3}\\
 && \langle \textsf{b}_1, u\rangle E_1+\langle \textsf{b}_2, u\rangle E_2 = 0, \label{eq4}\\
 && \langle f_1(u), v\rangle E_1+\langle f_2(u), v\rangle E_2 = 0, \label{eq5}\\
 && \langle \textsf{a}_2-\textsf{b}_1, u\rangle B_1 + \langle f_1(E_1), u\rangle - \langle \textsf{a}_1, E(u)-cu\rangle = 0, \label{eq6}\\
 && \langle \textsf{b}_1-\textsf{a}_2, u\rangle B_1 + \langle f_2(E_2), u\rangle - \langle \textsf{b}_2, E(u)-cu\rangle = 0, \label{eq7}\\
 && (B_2-A_1)\langle \textsf{a}_2, u\rangle +B_1 \langle \textsf{a}_1-\textsf{b}_2, u\rangle + \langle f_2(E_1), u\rangle - \langle \textsf{a}_2, E(u)-cu \rangle = 0, \label{eq8}\\
 && (B_2-A_1)\langle \textsf{b}_1, u\rangle +B_1 \langle \textsf{a}_1-\textsf{b}_2, u\rangle - \langle f_1(E_2), u\rangle + \langle \textsf{b}_1, E(u)-cu \rangle = 0, \label{eq9}\\
 && \langle f_1(u), v\rangle (A_1-c)+ \langle f_2(u), v\rangle B_1 +\langle \textsf{a}_1, v\rangle A(u) + \langle \textsf{b}_1, v\rangle B(u) - \langle f_1(E(u)-cu), v\rangle \nonumber\\
 && \hspace{5cm} - \langle \textsf{a}_1, u\rangle A(v) - \langle \textsf{b}_1, u\rangle B(v) + \langle f_1(E(v)-cv), u\rangle = 0, \label{eq10}\\
 && \langle f_2(u), v\rangle (B_2-c)+ \langle f_1(u), v\rangle B_1 +\langle \textsf{a}_2, v\rangle A(u) + \langle \textsf{b}_2, v\rangle B(u) - \langle f_2(E(u)-cu), v\rangle \nonumber\\
 && \hspace{5cm} - \langle \textsf{a}_2, u\rangle A(v) - \langle \textsf{b}_2, u\rangle B(v) + \langle f_2(E(v)-cv), u\rangle = 0.\label{eq11}
\end{eqnarray}
\end{theorem}

\begin{proof}
We need the Ricci transformation $\textsf{Ric}$ to use the Ricci soliton equation. But computing the Ricci transformation requires the Ricci tensor of $(G,\langle . , .\rangle)$. Suppose that $\textsf{ric}$ denotes the Ricci tensor of $(G,\langle . , .\rangle)$. We use the following formula (see \cite{Djiadeu-Boucetta-Wouafo}) to compute the Ricci tensor $\textsf{ric}$.
\begin{equation}\label{ricci tensor equation}
    \textsf{ric}(x,y)=-\frac{1}{2}\operatorname{tr}(\textsf{ad}_x\circ\textsf{ad}_y)-\frac{1}{2}\operatorname{tr}(\textsf{ad}_x\circ\textsf{ad}^{\ast}_y)-\frac{1}{4}\operatorname{tr}(\textsf{J}_x\circ\textsf{J}_y)
    -\frac{1}{2}\langle\textsf{ad}_Hx,y\rangle-\frac{1}{2}\langle\textsf{ad}_Hy,x\rangle,
\end{equation}
where $\langle H,x\rangle=\operatorname{tr}(\textsf{ad}_x)$ and $\textsf{J}_xy=\textsf{ad}^{\ast}_yx$, for all $x, y\in\frak{g}$. With some computations, for any $u, v\in\Gamma$, we have the following relations:
\begin{eqnarray*}
    && \operatorname{tr}(\textsf{ad}_{\textsf{e}_1}\circ\textsf{ad}_{\textsf{e}_1})=\operatorname{tr}(\textsf{ad}_{\textsf{e}_1}\circ\textsf{ad}_{\textsf{e}_2})=\operatorname{tr}(\textsf{ad}_{\textsf{e}_2}\circ\textsf{ad}_{\textsf{e}_1})
    =\operatorname{tr}(\textsf{ad}_{\textsf{e}_2}\circ\textsf{ad}_{\textsf{e}_2})=0,\\
    && \operatorname{tr}(\textsf{ad}_{\textsf{e}_1}\circ\textsf{ad}_{u})=\operatorname{tr}(\textsf{ad}_{u}\circ\textsf{ad}_{\textsf{e}_1})=\operatorname{tr}(\textsf{ad}_{\textsf{e}_2}\circ\textsf{ad}_{u})
    =\operatorname{tr}(\textsf{ad}_{u}\circ\textsf{ad}_{\textsf{e}_2})=0,\\
    && \operatorname{tr}(\textsf{ad}_{u}\circ\textsf{ad}_{v})=\langle\textsf{a}_1,u\rangle\langle\textsf{a}_1,v\rangle+\langle\textsf{b}_1,u\rangle\langle\textsf{a}_2,v\rangle
    +\langle\textsf{b}_1,v\rangle\langle\textsf{a}_2,u\rangle+\langle\textsf{b}_2,u\rangle\langle\textsf{b}_2,v\rangle,\\
    && \operatorname{tr}(\textsf{ad}_{\textsf{e}_1}\circ\textsf{ad}^\ast_{\textsf{e}_1})=\|\textsf{a}_1\|^2+\|\textsf{a}_2\|^2, \ \ \ \ \ \ \ \operatorname{tr}(\textsf{ad}_{\textsf{e}_1}\circ\textsf{ad}^\ast_{\textsf{e}_2})=\operatorname{tr}(\textsf{ad}_{\textsf{e}_2}\circ\textsf{ad}^\ast_{\textsf{e}_1})=\langle\textsf{a}_1,\textsf{b}_1\rangle+\langle\textsf{a}_2,\textsf{b}_2\rangle,\\
    && \operatorname{tr}(\textsf{ad}_{\textsf{e}_1}\circ\textsf{ad}^\ast_{u})=\operatorname{tr}(\textsf{ad}_{u}\circ\textsf{ad}^\ast_{\textsf{e}_1})=\langle u,f_1(\textsf{a}_1)+f_2(\textsf{a}_2)\rangle, \ \ \operatorname{tr}(\textsf{ad}_{\textsf{e}_2}\circ\textsf{ad}^\ast_{\textsf{e}_2})=\|\textsf{b}_1\|^2+\|\textsf{b}_2\|^2,\\
    && \operatorname{tr}(\textsf{ad}_{\textsf{e}_2}\circ\textsf{ad}^\ast_{u})=\operatorname{tr}(\textsf{ad}_{u}\circ\textsf{ad}^\ast_{\textsf{e}_2})=\langle u,f_1(\textsf{b}_1)+f_2(\textsf{b}_2)\rangle,\\
    && \operatorname{tr}(\textsf{ad}_{u}\circ\textsf{ad}^\ast_{v})=\langle\textsf{a}_1,u\rangle\langle\textsf{a}_1,v\rangle+\langle\textsf{b}_1,u\rangle\langle\textsf{b}_1,v\rangle
    +\langle\textsf{a}_2,u\rangle\langle\textsf{a}_2,v\rangle+\langle\textsf{b}_2,u\rangle\langle\textsf{b}_2,v\rangle\\
    && \hspace{5cm}+\langle f_1(u),f_1(v)\rangle+\langle f_2(u),f_2(v)\rangle,\\
    \end{eqnarray*}
    \begin{eqnarray*}
    && \operatorname{tr}(\textsf{J}_{\textsf{e}_1}\circ\textsf{J}_{\textsf{e}_1})=-2(\|\textsf{a}_1\|^2+\|\textsf{b}_1\|^2)+\operatorname{tr}(f^2_1), \\
    && \operatorname{tr}(\textsf{J}_{\textsf{e}_1}\circ\textsf{J}_{\textsf{e}_2})=\operatorname{tr}(\textsf{J}_{\textsf{e}_2}\circ\textsf{J}_{\textsf{e}_1})=-2(\langle \textsf{a}_1,\textsf{a}_2\rangle+\langle \textsf{b}_1,\textsf{b}_2\rangle)+\operatorname{tr}(f_1\circ f_2),\\
    && \operatorname{tr}(\textsf{J}_{\textsf{e}_2}\circ\textsf{J}_{\textsf{e}_2})=-2(\|\textsf{a}_2\|^2+\|\textsf{b}_2\|^2)+\operatorname{tr}(f^2_2),\\
    && \operatorname{tr}(\textsf{J}_{\textsf{e}_1}\circ\textsf{J}_{u})=\operatorname{tr}(\textsf{J}_{\textsf{e}_2}\circ\textsf{J}_{u})=\operatorname{tr}(\textsf{J}_{u}\circ\textsf{J}_{\textsf{e}_1})
    =\operatorname{tr}(\textsf{J}_{u}\circ\textsf{J}_{\textsf{e}_2})=\operatorname{tr}(\textsf{J}_{u}\circ\textsf{J}_{v})=0,\\
    && H=\textsf{a}_1+\textsf{b}_2,\\
    && \textsf{ad}_{H}\textsf{e}_1=(\|\textsf{a}_1\|^2+\langle \textsf{a}_1,\textsf{b}_2\rangle)\textsf{e}_1+\langle \textsf{a}_2,\textsf{a}_1+\textsf{b}_2\rangle\textsf{e}_2,\\
    && \textsf{ad}_{H}\textsf{e}_2=(\|\textsf{b}_2\|^2+\langle \textsf{a}_1,\textsf{b}_2\rangle)\textsf{e}_2+\langle \textsf{b}_1,\textsf{a}_1+\textsf{b}_2\rangle\textsf{e}_1,\\
    && \textsf{ad}_{H}u=\langle f_1(\textsf{a}_1+\textsf{b}_2), u \rangle\textsf{e}_1+\langle f_2(\textsf{a}_1+\textsf{b}_2), u \rangle\textsf{e}_2.
\end{eqnarray*}
Now the above relations together with the formula \eqref{ricci tensor equation} give the following formulas for the Ricci tensor $\textsf{ric}$.
\begin{eqnarray*}
  \textsf{ric}(\textsf{e}_1,\textsf{e}_1) &=& \frac{1}{2}(\|\textsf{b}_1\|^2-\|\textsf{a}_2\|^2)-\frac{1}{4}\operatorname{tr}(f^2_1)-\|\textsf{a}_1\|^2-\langle \textsf{a}_1,\textsf{b}_2\rangle, \\
  \textsf{ric}(\textsf{e}_1,\textsf{e}_2) &=& -(\langle \textsf{a}_1,\textsf{b}_1\rangle+\langle \textsf{a}_2,\textsf{b}_2\rangle+\frac{1}{4}\operatorname{tr}(f_1\circ f_2)), \\
  \textsf{ric}(\textsf{e}_1,u) &=&  \langle u, -f_1(\textsf{a}_1)-\frac{1}{2}(f_2(\textsf{a}_2)+f_1(\textsf{b}_2))\rangle,\\
  \textsf{ric}(\textsf{e}_2,\textsf{e}_2) &=& \frac{1}{2}(\|\textsf{a}_2\|^2-\|\textsf{b}_1\|^2)-\frac{1}{4}\operatorname{tr}(f^2_2)-\|\textsf{b}_2\|^2-\langle \textsf{a}_1,\textsf{b}_2\rangle, \\
  \textsf{ric}(\textsf{e}_2,u) &=&  \langle u, -f_2(\textsf{b}_2)-\frac{1}{2}(f_1(\textsf{b}_1)+f_2(\textsf{a}_1))\rangle,\\
  \textsf{ric}(u,v) &=& -\langle \textsf{a}_1, u \rangle\langle \textsf{a}_1, v \rangle-\langle \textsf{b}_2, u \rangle\langle \textsf{b}_2, v \rangle-\frac{1}{2}\langle u,\textsf{b}_1+\textsf{a}_2 \rangle\langle v,\textsf{b}_1+\textsf{a}_2 \rangle+\frac{1}{2}\langle u,(f^2_1+f^2_2)(v) \rangle.
\end{eqnarray*}
We know that the relation between Ricci transformation $\textsf{Ric}$ and the Ricci tensor $\textsf{ric}$ is given by the equation $\langle \textsf{Ric}(x), y \rangle=\textsf{ric}(x,y)$, for all $x,y\in\frak{g}$. Hence, we obtain the Ricci transformation as follows:
\begin{eqnarray*}
  \textsf{Ric}(\textsf{e}_1) &=& A_1\textsf{e}_1+B_1\textsf{e}_2+E_1, \\
  \textsf{Ric}(\textsf{e}_2) &=& A_2\textsf{e}_1+B_2\textsf{e}_2+E_2, \\
  \textsf{Ric}(u) &=& A(u)\textsf{e}_1+B(u)\textsf{e}_2+E(u),
\end{eqnarray*}
where $u\in\Gamma$. The above equations, together with $\textsf{Ric}(x)=cx+D(x)$ and the fact that $D$ is a derivation, complete the proof.
\end{proof}
\begin{cor}\label{corollary two-step nilpotent two-dim commutator}
Let $G\in\textsf{Lie(n,2)}$ be a $2$-step nilpotent Lie group. Then, $(G,\langle . , .\rangle)$ is a Ricci soliton if and only if the following two equations hold for any $u\in\Gamma$:
\begin{eqnarray*}
  cf_1(u) &=&  (\frac{1}{4}\operatorname{tr}(f^2_1))f_1(u)+(\frac{1}{4}\operatorname{tr}(f_1\circ f_2))f_2(u)+\frac{1}{2}f_1((f^2_1+f^2_2)(u))+\frac{1}{2}(f^2_1+f^2_2)(f_1(u)),\\
  cf_2(u) &=& (\frac{1}{4}\operatorname{tr}(f^2_2))f_2(u)+(\frac{1}{4}\operatorname{tr}(f_1\circ f_2))f_1(u)+\frac{1}{2}f_2((f^2_1+f^2_2)(u))+\frac{1}{2}(f^2_1+f^2_2)(f_2(u)).
\end{eqnarray*}
\end{cor}
\begin{proof}
If $G$ is $2$-step nilpotent, then $\textsf{a}_1=\textsf{a}_2=\textsf{b}_1=\textsf{b}_2=0$. So we have:
\begin{equation*}
    A_1=-\frac{1}{4}\operatorname{tr}(f^2_1), A_2=B_1=-\frac{1}{4}\operatorname{tr}(f_1\circ f_2), B_2=-\frac{1}{4}\operatorname{tr}(f^2_2), E_1=E_2=0.
\end{equation*}
Also, for any $u\in\Gamma$ we have:
\begin{equation*}
    A(u)=B(u)=0, \ \ \ \ E(u)=\frac{1}{2}(f^2_1+f^2_2)(u).
\end{equation*}
In this case, the first nine equations of the previous theorem hold automatically, so it suffices to check the last two equations, which reduce to the two equations in this corollary.
\end{proof}
\begin{example}\label{example 1 two-step nilpotent}
Let $N$ be the two-step nilpotent Lie group whose Lie algebra $\frak{n}:=\textsf{span}{e_1,\dots,e_n}$, $(n\geq4)$, has the following non-zero Lie bracket:
\begin{equation*}
    [e_3,e_4]=e_1+e_2.
\end{equation*}
Suppose that $\langle \cdot,\cdot \rangle$ is an inner product on $\frak{n}$ such that ${e_1,\dots,e_n}$ is an orthonormal basis. We can easily see that $\textsf{a}_1=\textsf{a}_2=\textsf{b}_1=\textsf{b}_2=0$, $\Gamma=\textsf{span}{e_3,\dots,e_n}$, and $f_1=f_2$ have the following representation with respect to the basis $\mathcal{B}={e_3,\dots,e_n}$:
\begin{equation*}
    f_1=f_2= \left(
  \begin{array}{ccccc}
    0 & -1 & 0 & \dots & 0 \\
    1 & 0 & 0 & \dots & 0 \\
    0 & 0 & 0 & \dots & 0 \\
    \vdots & \vdots & \vdots &  & \vdots \\
    0 & 0 & 0 & \dots & 0 \\
  \end{array}
\right).
\end{equation*}
If we define $f:=f_1=f_2$, then the two relations given in Corollary \ref{corollary two-step nilpotent two-dim commutator} reduce to the following simple equation:
\begin{equation*}
    cf(u)=(\frac{1}{2}\operatorname{tr}f^2)f(u)+2f^3(u).
\end{equation*}
It follows that $f^3=-f$ and $f^2$ has the following representation with respect to the basis $\mathcal{B}$:
\begin{equation*}
    f^2= \left(
  \begin{array}{ccccc}
    -1 & 0 & 0 & \dots & 0 \\
    0 & -1 & 0 & \dots & 0 \\
    0 & 0 & 0 & \dots & 0 \\
    \vdots & \vdots & \vdots &  & \vdots \\
    0 & 0 &  0 & \dots & 0 \\
  \end{array}
\right).
\end{equation*}
Thus, $(\frak{n},\langle\cdot,\cdot\rangle)$ is a two-step nilpotent Ricci nilsoliton with $c=-3$.
\end{example}

\begin{example}\label{example 2 two-step nilpotent}
Suppose that $N$ is an $n$-dimensional two-step nilpotent Lie group, where $n$ is an even natural number greater than or equal to $4$, with the following non-zero Lie brackets for its Lie algebra $\frak{n}=\textsf{span}{e_1,\dots,e_n}$:
\begin{equation*}
    [e_{2i+1},e_{2i+2}]=e_1+e_2, \ \ \ \ i=1,\dots, \frac{n-2}{2}.
\end{equation*}
Let $\langle\cdot,\cdot\rangle$ be the inner product defined on $\frak{n}$ such that ${e_1,\dots,e_n}$ is an orthonormal basis. Then $\textsf{a}_1=\textsf{a}_2=\textsf{b}_1=\textsf{b}_2=0$, $\Gamma=\textsf{span}{e_3,\dots,e_n}$, and $f_1=f_2$ have the following representation with respect to the basis $\mathcal{B}={e_3,\dots,e_n}$:
\begin{equation*}
    f_1=f_2= \left(
  \begin{array}{ccccccc}
    0 & -1 & 0 & 0 & \dots & 0 & 0 \\
    1 & 0 & 0 & 0 & \dots & 0 & 0 \\
    0 & 0 & 0 & -1 & \dots & 0 & 0 \\
    0 & 0 & 1 & 0 & \dots & 0 & 0 \\
    \vdots & \vdots & \vdots & \vdots &  & \vdots & \vdots \\
    0 & 0 & 0 & 0 &  \dots & 0 & -1 \\
    0 & 0 & 0 & 0 &  \dots & 1 & 0 \\
  \end{array}
\right).
\end{equation*}
We can easily see that $f:=f_1=f_2$ is an almost complex structure on $\Gamma$. If we extend the almost complex structure $f$ to an almost complex structure $J$ on $\frak{n}$ by $Je_1=e_2$ and $Je_2=-e_1$, then using Theorems 2.1 and 3.1 of \cite{Salimi Moghaddam}, $(J,\langle\cdot,\cdot\rangle)$ is a non-K\"ahler Hermitian structure on $N$. Now we can see that the equations of Corollary \ref{corollary two-step nilpotent two-dim commutator} reduce to the equation $cf=\left(\frac{1}{2}\operatorname{tr}f^2\right)f+2f^3$. In this case we have $f^2=-\operatorname{Id}\Gamma$ and $f^3=-f$. Thus, $(N,J,\langle\cdot,\cdot\rangle)$ is a non-K\"ahler Hermitian Ricci nilsoliton with $c=-\frac{n+2}{2}$.
\end{example}

\section{\textbf{Indecomposable Lie groups}}\label{Indecomposable Lie groups}
In this section, we apply Theorem \ref{Classi Two-dim commu Ricci solitons} to simply connected Lie groups whose Lie algebras are indecomposable $G\in\textsf{Lie(n,2)}$ with $n\geq 5$.
We use the results of \cite{Le-Nguyen-Vo}, which classifies all such indecomposable Lie groups $G$ into three types. We examine the results of section \ref{2-dimensional} for these three classes.

\begin{example}\label{example 1}
Suppose that $G\in\textsf{Lie(n,2)}$ is an indecomposable Lie group with $n=5+2k$, such that $k\geq 0$. Theorem 10 of \cite{Le-Nguyen-Vo} shows that there exists a basis $\{X_1, X_2,\dots, X_n\}$ for the Lie algebra $\frak{g}$ such that the non-zero Lie brackets are as follows:
\begin{equation}\label{Lie brackets of the case n=5+2k}
    [X_3,X_4]=X_1,\ \ \ \ \ [X_3,X_1]=[X_4,X_5]=\dots=[X_{4+2k},X_{5+2k}]=X_2.
\end{equation}
Let $\langle \cdot , \cdot \rangle$ be the left-invariant Riemannian metric on $G$ for which $\{X_1, X_2,\dots, X_n\}$  is an orthonormal basis. With respect to the induced inner product on $\frak{g}$, we can consider $\textsf{e}_1=X_1$ and $\textsf{e}_2=X_2$. So we have $\Gamma=\textsf{span}\{X_3, X_4,\dots, X_n\}$. Using the Lie algebra structure equations given in \eqref{Lie brackets of the case n=5+2k}, we have $\textsf{a}_2=X_3$ and $\textsf{a}_1=\textsf{b}_1=\textsf{b}_2=0$. Also, the matrix representation of the linear maps $f_1$ and $f_2$, with respect to the basis $\{X_3, X_4,\dots, X_n\}$, are of the following forms:
\begin{equation}\label{f_1 and f_2 matrix case n=5+2k}
    f_1= \left(
  \begin{array}{ccccc}
    0 & -1 & 0 & \dots & 0 \\
    1 & 0 & 0 & \dots & 0 \\
    0 & 0 & 0 & \dots & 0 \\
    \vdots & \vdots & \vdots &  & \vdots \\
    0 & 0 & 0 & \dots & 0 \\
  \end{array}
\right), \ \ \ \ \
f_2= \left(
  \begin{array}{cccccc}
    0 & 0 & 0 & \dots & 0 & 0 \\
    0 & 0 & -1 & \dots & 0 & 0 \\
    0 & 1 & 0 & \dots & 0 & 0 \\
    \vdots & \vdots & \vdots &  & \vdots &\vdots \\
    0 & 0 & 0 & \dots & 0 & -1 \\
    0 & 0 & 0 & \dots & 1 & 0 \\
  \end{array}
\right).
\end{equation}
Using Proposition \ref{Levi-Civita and Sectional}, the Levi-Civita connection is given in Table \ref{table 2}.
\begin{table}[h]
\caption{Levi-Civita connection of an indecomposable Lie group $G\in\textsf{Lie(5+2k,2)}$}
\label{table 2}
\centering
\begin{tabular}{|c|c|c|c|c|c|c|c|c|}
  \hline
  $\nabla$ & $\textsf{e}_1$ & $\textsf{e}_2$ & $X_3$ & $X_4$ & $X_5$ & $\dots$ & $X_{4+2k}$ & $X_{5+2k}$ \\
  \hline
  $\textsf{e}_1$ & $0$ & $\frac{1}{2}X_3$ & $-\frac{1}{2}(\textsf{e}_2+X_4)$ & $\frac{1}{2}X_3$ & $0$ & $\dots$ & $0$ & $0$ \\
  \hline
  $\textsf{e}_2$ & $\frac{1}{2}X_3$ & $0$ & $-\frac{1}{2}\textsf{e}_1$ & $-\frac{1}{2}X_5$ & $\frac{1}{2}X_4$ & $\dots$ & $-\frac{1}{2}X_{5+2k}$ & $\frac{1}{2}X_{4+2k}$ \\
  \hline
  $X_3$ & $\frac{1}{2}(\textsf{e}_2-X_4)$ & $-\frac{1}{2}\textsf{e}_1$ & $0$ & $\frac{1}{2}\textsf{e}_1$ & $0$ & $\dots$ & $0$ & $0$ \\
  \hline
  $X_4$ & $\frac{1}{2}X_3$ & $-\frac{1}{2}X_5$ & $-\frac{1}{2}\textsf{e}_1$ & $0$ & $\frac{1}{2}\textsf{e}_2$ & $\dots$ & $0$ & $0$ \\
  \hline
  $X_5$ & $0$ & $\frac{1}{2}X_4$ & $0$ & $-\frac{1}{2}\textsf{e}_2$ & $0$ & $\dots$ & $0$ & $0$ \\
  \hline
  $\vdots$ & $ $ & $ $ & $ $ & $ $ & $ $ & $ $ & $\vdots$ & $\vdots$ \\
  \hline
  $X_{4+2k}$ & $0$ & $-\frac{1}{2}X_{5+2k}$ & $0$ & $0$ & $0$ & $\dots$ & $0$ & $\frac{1}{2}\textsf{e}_2$ \\
  \hline
  $X_{5+2k}$ & $0$ & $\frac{1}{2}X_{4+2k}$ & $0$ & $0$ & $0$ & $\dots$ & $-\frac{1}{2}\textsf{e}_2$ & $0$ \\
  \hline
\end{tabular}
\end{table}

In this case, a direct computation shows that $\operatorname{tr}(f^2_1)=-2, \operatorname{tr}(f^2_2)=3-n$ and $\operatorname{tr}(f_1\circ f_2)=0$. Thus, the Ricci transformation of $(G,\langle \cdot , \cdot \rangle)$ is:
\begin{equation*}
    \textsf{Ric}(\textsf{e}_1) = 0, \ \ \textsf{Ric}(\textsf{e}_2) = \frac{n-1}{4}\textsf{e}_2, \ \ \textsf{Ric}(u) = -\frac{\lambda_3}{2}X_3+\frac{1}{2}\left(
                                                                                                                                                          \begin{array}{ccccc}
                                                                                                                                                            -1 & 0 & 0 & \dots & 0 \\
                                                                                                                                                            0 & -2 & 0 & \dots & 0 \\
                                                                                                                                                            0 & 0 & -1 & \dots & 0 \\
                                                                                                                                                            \vdots & \vdots & \vdots &  & \vdots \\
                                                                                                                                                            0 & 0 & 0 & \dots & -1 \\
                                                                                                                                                          \end{array}
                                                                                                                                                        \right)u,
\end{equation*}
where $u=\sum_{i=3}^n\lambda_iX_i$ is an arbitrary element of $\Gamma$, and the matrix is considered with respect to the basis $\{X_3, X_4,\dots, X_n\}$.\\
The equation \eqref{eq8} of Theorem \ref{Classi Two-dim commu Ricci solitons} shows that $c=-\frac{n+3}{4}$.
On the other hand, taking $u=X_4$ and $v=X_5$ in \eqref{eq11} gives $0=\frac{1}{2}$, which is impossible.
\end{example}

If $n=6+2k$, then there are two different Lie algebra structures on $\frak{g}$. We discuss these Lie groups in the following two examples.

\begin{example}\label{example 2}
Let $G\in\textsf{Lie(n,2)}$ be an indecomposable Lie group such that $n=6+2k$, where $k\geq 0$. Then there are two Lie algebra structures. We study the first one in this example. Using \cite{Le-Nguyen-Vo}, there is a basis $\{X_1, X_2,\dots, X_n\}$ for the Lie algebra $\frak{g}$ such that:
\begin{equation}\label{Lie brackets of the case n=6+2k case 1}
    [X_3,X_1]=X_1,\ \ \ \ \ [X_3,X_4]=[X_5,X_6]=\dots=[X_{5+2k},X_{6+2k}]=X_2.
\end{equation}
Assume that $\langle \cdot , \cdot \rangle$ is a left-invariant Riemannian metric on $G$ such that $\{X_1, X_2,\dots, X_n\}$ is orthonormal. Similar to the previous example we set $\textsf{e}_1=X_1$ and $\textsf{e}_2=X_2$. Therefore  $\Gamma=\textsf{span}\{X_3, X_4,\dots, X_n\}$. Using \eqref{Lie brackets of the case n=6+2k case 1}, we have $\textsf{a}_1=X_3$ and $\textsf{a}_2=\textsf{b}_1=\textsf{b}_2=0$. We see that $f_1=0$ and the matrix representation of $f_2$ is as follows:
\begin{equation}\label{f_2 matrix case n=6+2k case 1}
   f_2= \left(
  \begin{array}{ccccccc}
    0 & -1 & 0 & 0 &\dots & 0 & 0 \\
    1 & 0 & 0 & 0 & \dots & 0 & 0 \\
    0 & 0 & 0 & -1 & \dots & 0 & 0 \\
    0 & 0 & 1 & 0 & \dots & 0 & 0 \\
    \vdots & \vdots & \vdots & \vdots &  & \vdots &\vdots \\
    0 & 0 & 0 & 0 & \dots & 0 & -1 \\
    0 & 0 & 0 & 0 & \dots & 1 & 0 \\
  \end{array}
\right).
\end{equation}
The Levi-Civita connection is given in Table \ref{table 3}.
\begin{table}[h]
\caption{Levi-Civita connection of an indecomposable Lie group $G\in\textsf{Lie(6+2k,2)}$ (Case 1)}
\label{table 3}
\centering
\begin{tabular}{|c|c|c|c|c|c|c|c|c|c|}
  \hline
  $\nabla$ & $\textsf{e}_1$ & $\textsf{e}_2$ & $X_3$ & $X_4$ & $X_5$ & $X_6$ & $\dots$ & $X_{5+2k}$ & $X_{6+2k}$ \\
  \hline
  $\textsf{e}_1$ & $X_3$ & $0$ & $-\textsf{e}_1$ & $0$ & $0$ & $0$ & $\dots$ & $0$ & $0$ \\
  \hline
  $\textsf{e}_2$ & $0$ & $0$ & $-\frac{1}{2}X_4$ & $\frac{1}{2}X_3$ & $-\frac{1}{2}X_6$ & $\frac{1}{2}X_5$ & $\dots$ & $-\frac{1}{2}X_{6+2k}$ & $\frac{1}{2}X_{5+2k}$ \\
  \hline
  $X_3$ & $0$ & $-\frac{1}{2}X_4$ & $0$ & $\frac{1}{2}\textsf{e}_2$ & $0$ & $0$ & $\dots$ & $0$ & $0$ \\
  \hline
  $X_4$ & $0$ & $\frac{1}{2}X_3$ & $-\frac{1}{2}\textsf{e}_2$ & $0$ & $0$ & $0$ & $\dots$ & $0$ & $0$ \\
  \hline
  $X_5$ & $0$ & $-\frac{1}{2}X_6$ & $0$ & $0$ & $0$ & $\frac{1}{2}\textsf{e}_2$ & $\dots$ & $0$ & $0$ \\
  \hline
  $X_6$ & $0$ & $\frac{1}{2}X_5$ & $0$ & $0$ & $-\frac{1}{2}\textsf{e}_2$ & $0$ & $\dots$ & $0$ & $0$ \\
  \hline
  $\vdots$ & $ $ & $ $ & $ $ & $ $ & $ $ & $ $ & $ $ & $\vdots$ & $\vdots$ \\
  \hline
  $X_{5+2k}$ & $0$ & $-\frac{1}{2}X_{6+2k}$ & $0$ & $0$ & $0$ & $0$ & $\dots$ & $0$ & $\frac{1}{2}\textsf{e}_2$ \\
  \hline
  $X_{6+2k}$ & $0$ & $\frac{1}{2}X_{5+2k}$ & $0$ & $0$ & $0$ & $0$ & $\dots$ & $-\frac{1}{2}\textsf{e}_2$ & $0$ \\
  \hline
\end{tabular}
\end{table}

We have $\operatorname{tr}(f^2_1)=0, \operatorname{tr}(f^2_2)=-(n-2)$, and $\operatorname{tr}(f_1\circ f_2)=0$. Therefore the Ricci transformation is as follows:
\begin{equation*}
    \textsf{Ric}(\textsf{e}_1) = -\textsf{e}_1, \ \ \textsf{Ric}(\textsf{e}_2) = \frac{n-2}{4}\textsf{e}_2-\frac{1}{2}X_4, \ \ \textsf{Ric}(u) = -\frac{1}{2}\langle u,X_4\rangle\textsf{e}_2-\langle u,X_3\rangle X_3-\frac{1}{2}u,
\end{equation*}
for all $u\in\Gamma$.\\
Using \eqref{eq6} of Theorem \ref{Classi Two-dim commu Ricci solitons} with $u=X_3$, we obtain
$c=-1$. However, substituting $c=-1$, $u=X_3$ and $v=X_4$ into equation \eqref{eq11} gives $n=-2$ which is impossible.
\end{example}

In the following example, we study the second type of $n$-dimensional indecomposable Lie groups $G\in\textsf{Lie(n,2)}$, where $n=6+2k$.

\begin{example}\label{example 3}
Assume that $G\in\textsf{Lie(n,2)}$ is the second type of $n$-dimensional indecomposable Lie group, classified in \cite{Le-Nguyen-Vo}, with $n=6+2k$, $k\geq 0$. Then there exists a basis $\{X_1, X_2,\dots, X_n\}$ for the Lie algebra $\frak{g}$ such that:
\begin{equation}\label{Lie brackets of the case n=6+2k case 2}
    [X_3,X_4]=X_1,\ \ \ \ \ [X_3,X_1]=[X_5,X_6]=\dots=[X_{5+2k},X_{6+2k}]=X_2.
\end{equation}
Let $\langle \cdot , \cdot \rangle$ be the left-invariant Riemannian metric on $G$ for which $\{X_1, X_2,\dots, X_n\}$ is orthonormal.
In this case, $\textsf{e}_1=X_1$, $\textsf{e}_2=X_2$, and $\Gamma=\textsf{span}\{X_3, X_4,\dots, X_n\}$.
Using \eqref{Lie brackets of the case n=6+2k case 2}, we have $\textsf{a}_2=X_3$ and $\textsf{a}_1=\textsf{b}_1=\textsf{b}_2=0$.
The matrix representation of $f_1$ is the same as in Example \ref{example 1}, and the matrix representation of $f_2$, with respect to the basis $\{X_3, X_4,\dots, X_n\}$, is of the following form:
\begin{equation}\label{f_2 matrix case n=6+2k case 2}
   f_2= \left(
  \begin{array}{ccccccc}
    0 & 0 & 0 & 0 & \dots & 0 & 0 \\
    0 & 0 & 0 & 0 & \dots & 0 & 0 \\
    0 & 0 & 0 & -1 & \dots & 0 & 0 \\
    0 & 0 & 1 & 0 & \dots & 0 & 0 \\
    \vdots & \vdots & \vdots & \vdots &  & \vdots &\vdots \\
    0 & 0 & 0 & 0 & \dots & 0 & -1 \\
    0 & 0 & 0 & 0 & \dots & 1 & 0 \\
  \end{array}
\right).
\end{equation}
The Levi-Civita connection of this Lie group is as shown in table \ref{table 4}.
\begin{table}[h]
\caption{Levi-Civita connection of an indecomposable Lie group $G\in\textsf{Lie(6+2k,2)}$ (Case 2)}
\label{table 4}
\centering
\begin{tabular}{|c|c|c|c|c|c|c|c|c|c|}
  \hline
  $\nabla$ & $\textsf{e}_1$ & $\textsf{e}_2$ & $X_3$ & $X_4$ & $X_5$ & $X_6$ & $\dots$ & $X_{5+2k}$ & $X_{6+2k}$ \\
  \hline
  $\textsf{e}_1$ & $0$ & $\frac{1}{2}X_3$ & $-\frac{1}{2}(\textsf{e}_2+X_4)$ & $\frac{1}{2}X_3$ & $0$ & $0$ & $\dots$ & $0$ & $0$ \\
  \hline
  $\textsf{e}_2$ & $\frac{1}{2}X_3$ & $0$ & $-\frac{1}{2}\textsf{e}_1$ & $0$ & $-\frac{1}{2}X_6$ & $\frac{1}{2}X_5$ & $\dots$ & $-\frac{1}{2}X_{6+2k}$ & $\frac{1}{2}X_{5+2k}$ \\
  \hline
  $X_3$ & $\frac{1}{2}(\textsf{e}_2-X_4)$ & $-\frac{1}{2}\textsf{e}_1$ & $0$ & $\frac{1}{2}\textsf{e}_1$ & $0$ & $0$ & $\dots$ & $0$ & $0$ \\
  \hline
  $X_4$ & $\frac{1}{2}X_3$ & $0$ & $-\frac{1}{2}\textsf{e}_1$ & $0$ & $0$ & $0$ & $\dots$ & $0$ & $0$ \\
  \hline
  $X_5$ & $0$ & $-\frac{1}{2}X_6$ & $0$ & $0$ & $0$ & $\frac{1}{2}\textsf{e}_2$ & $\dots$ & $0$ & $0$ \\
  \hline
  $X_6$ & $0$ & $\frac{1}{2}X_5$ & $0$ & $0$ & $-\frac{1}{2}\textsf{e}_2$ & $0$ & $\dots$ & $0$ & $0$ \\
  \hline
  $\vdots$ & $ $ & $ $ & $ $ & $ $ & $ $ & $ $ & $ $ & $\vdots$ & $\vdots$ \\
  \hline
  $X_{5+2k}$ & $0$ & $-\frac{1}{2}X_{6+2k}$ & $0$ & $0$ & $0$ & $0$ & $\dots$ & $0$ & $\frac{1}{2}\textsf{e}_2$ \\
  \hline
  $X_{6+2k}$ & $0$ & $\frac{1}{2}X_{5+2k}$ & $0$ & $0$ & $0$ & $0$ & $\dots$ & $-\frac{1}{2}\textsf{e}_2$ & $0$ \\
  \hline
\end{tabular}
\end{table}

Using the representations of $f_1$ and $f_2$ we have $\operatorname{tr}(f^2_1)=-2, \operatorname{tr}(f^2_2)=4-n$ and $\operatorname{tr}(f_1\circ f_2)=0$. So for the Ricci transformation, we have:
\begin{equation*}
    \textsf{Ric}(\textsf{e}_1) = 0, \ \ \textsf{Ric}(\textsf{e}_2) = \frac{k+1}{2}\textsf{e}_2, \ \ \textsf{Ric}(u) = -\frac{1}{2}(\langle u,X_3\rangle X_3+u),
\end{equation*}
for all $u\in\Gamma$.\\
Now, \eqref{eq8} of Theorem \ref{Classi Two-dim commu Ricci solitons}, with $u=X_3$, gives $c=-\frac{k+3}{2}$. If we substitute $c=-\frac{k+3}{2}$ into equation \eqref{eq10} and set $u=X_3$, $v=X_4$, we obtain $c=-\frac{3}{2}$.
On the other hand, using \eqref{eq8} we have $k=0$. So for $k\neq0$, $(G,\langle \cdot , \cdot \rangle)$ is not a Ricci soliton. But if we consider the case $k=0$ (or equivalently $n=6$) then we have:
\begin{equation*}
    \textsf{Ric}(\textsf{e}_1) = 0, \ \ \textsf{Ric}(\textsf{e}_2) = \frac{1}{2}\textsf{e}_2, \ \ \textsf{Ric}(X_3) = -X_3, \ \ \textsf{Ric}(X_i) = -\frac{1}{2}X_i, \ \ i=4,5,6.
\end{equation*}

For $k=0$ (equivalently, $n=6$), let $c=-\frac32$ and define
$D\in\operatorname{Der}(\mathfrak g)$ by
\[
D(\textsf{e}_1)=\frac{3}{2}\textsf{e}_1,\qquad
D(\textsf{e}_2)=2\textsf{e}_2,\qquad
D(X_3)=\frac{1}{2}X_3,\qquad
D(X_i)=X_i,\quad i=4,5,6.
\]
Then $\operatorname{Ric}=c\,\operatorname{Id}+D$. Hence, $(G,\langle\cdot,\cdot\rangle)$ is a Ricci soliton.
\end{example}

\section*{Declarations}

{\textbf{Conflict of interest:} There is no conflict of interest.\\

{\textbf{Ethical approval:} Not applicable.\\

{\textbf{Competing interests:} The author has no relevant financial or non-financial interests to disclose.\\

{\textbf{Authors contributions:}} Not applicable.\\

{\textbf{Funding:}} The author received no financial support. \\

{\textbf{Availability of data and materials:}} Data sharing not applicable to this article as no datasets were generated or analysed during the current study.\\


\bibliographystyle{amsplain}

\begin{thebibliography}{99}

\bibitem{Asgari-Salimi} F. Asgari and H. R. Salimi Moghaddam, \emph{Riemannian geometry of two families of tangent Lie groups}, Bull. Iran. Math. Soc., \textbf{44}(2018), 193--203.

\bibitem{Chow-Knopf} B. Chow and D. Knopf, \emph{The Ricci flow: an introduction}, AMS, Math. Surveys Monogr., \textbf{110}, (2004).

\bibitem{Djiadeu-Boucetta-Wouafo} M.B. Djiadeu Ngaha, M. Boucetta and J. Wouafo Kamga, \emph{The signature of the Ricci curvature of left-invariant Riemannian metrics on nilpotent Lie groups}, Differential Geom. Appl., \textbf{47}(2016), 26--42.

\bibitem{Gordon-Jablonski2024} C. Gordon and M. Jablonski, \emph{Infinitesimal maximal symmetry and Ricci soliton solvmanifolds}, Trans. Amer. Math. Soc., \textbf{377}(8)(2024), 5673--5704.

\bibitem{Hamilton} R. S. Hamilton, \emph{Three-manifolds with positive Ricci curvature}, J. Differential Geom., \textbf{17}(2)(1982), 255--306.

\bibitem{Jablonski1} M. Jablonski, \emph{Detecting orbits along subvarieties via the moment map}, M\"unster J. of Math., \textbf{3}(2010), 67--88.


\bibitem{Jablonski2} M. Jablonski, \emph{Homogeneous Ricci solitons are algebraic}, Geom. Topol., \textbf{18}(2014), 2477--2486.

\bibitem{Jablonski3} M. Jablonski, \emph{Homogeneous Einstein manifolds}, Rev. Un. Mat. Argentina, \textbf{64}(2)(2023), 461--485.

\bibitem{Kaiser} V. Kaiser, \emph{Ricci curvatures of left invariant metrics on Lie groups with one-dimensional commutator groups}, Note Mat., \textbf{14}(2)(1994), 209-215.

\bibitem{Lauret1} J. Lauret, \emph{Ricci soliton homogeneous nilmanifolds}, Math. Ann., \textbf{319}(2001), 715--733.

\bibitem{Lauret2} J. Lauret, \emph{Einstein solvmanifolds are standard}, Ann. Math., \textbf{172}(3)(2010), 1859--1877.

\bibitem{Le-Nguyen-Vo} V. A. Le, T. A. Nguyen, T. T. C. Nguyen, T. T. M. Nguyen and T. N. Vo , \emph{Applying matrix theory to classify real solvable Lie algebras having 2-dimensional derived ideals}, Linear Algebra Appl., \textbf{588}(2020), 282--303.

\bibitem{Mi} J. Milnor, \emph{Curvatures of left invariant metrics on Lie groups}, Adv. Math., \textbf{21} (1976), 293--329.

\bibitem{Nikolayevsky1} Y. Nikolayevsky, \emph{Einstein solvmanifolds and the Pre-Einstein Derivation}, Trans. Am. Math. Soc., \textbf{363}(8)(2011), 3935--3958.

\bibitem{Nikolayevsky2} Y. Nikolayevsky, \emph{Einstein solvmanifolds attached to two-step nilradicals}, Math. Z., \textbf{272} (2012), 675--695.

\bibitem{Salimi Moghaddam} H. R. Salimi Moghaddam, \emph{Hermitian and K\"ahler Structures of Two Types on Lie Groups with Two-Dimensional Derived Algebras},
https://doi.org/10.48550/arXiv.2602.13221

\end{thebibliography}

\end{document}